\documentclass[a4paper,draft]{amsproc}

\usepackage[alphabetic,y2k,initials]{amsrefs}
\usepackage{amsmath,amssymb,amsfonts}
\usepackage[dvips]{graphicx}

\usepackage{pst-all}
\usepackage{pstricks-add}
\usepackage{pst-slpe}
\usepackage{pst-solides3d}

\usepackage[all]{xy}
\usepackage{textcomp}

\theoremstyle{plain}
 
 \newtheorem{prop}{Proposition}[section]

\theoremstyle{definition}

\theoremstyle{remark}
 \newtheorem{rem}{Remark}[section]
 \numberwithin{equation}{section}

\renewcommand{\le}{\leqslant}

\title[Integrable lattices of hyperplanes related to billiards]{Integrable lattices of hyperplanes related to billiards within confocal quadrics}

 \subjclass[2000]{37J35, 14H70, 37D05, 51N25}
 \keywords{billiard, confocal quadrics, discrete integrable systems, double reflection nets}

\author[Radnovi\'c]{\bfseries Milena Radnovi\'c}

\address{Mathematical Institute SANU, Kneza Mihaila 36, Belgrade, Serbia
\newline\indent School of Mathematics and Statistics, University of Sydney, NSW, Australia}

\email{milena@mi.sanu.ac.rs, milena.radnovic@sydney.edu.au}

\thanks{The research which lead to this paper was partially
supported by the Serbian Ministry of Education and Science (Project
no.~174020: \emph{Geometry and Topology of Manifolds, Classical Mechanics and Integrable Dynamical Systems})
and by grant no.~FL120100094 from the Australian Research Council.
}

\date{}

\begin{document}

\begin{abstract}
We introduce a new discrete system that arises from ellipsoidal billiards and is closely related to the double reflection nets.
The system is defined on the lattice of a uniform honeycomb consisting of rectified hypercubes and cross polytopes.
In the $2$-dimensional case, the lattice is regular and it incorporates dynamics both in the original space and its dual.
In the $3$-dimensional case, the lattice consists of tetrahedra and cuboctahedra.
\end{abstract}

\maketitle
\tableofcontents

\section{Introduction}
Dicrete integrable systems represent a few decades old field of study that excites both mathematics and physics communities.
Mathematical areas of differential equations, dynamical systems, geometry, algebraic geometry are all intertwined in that field, which, on the other hand, has applications in various physics theories. 
For detailed expostions and references, see for example \cites{GrammaticosDIS, BS2008book, DuistermaatBOOK} and referenced therein.

This paper is focused to discrete systems arising from the dynamics of billiards within confocal quadrics.
A class of such systems, \emph{double reflection nets}, is introduced and studied in \cites{DragRadn2012jnmp,DragRadn2014bul,DragRadn2015umn}.
Double reflection nets are \emph{discrete line congruences}, that is maps that assigns a line to each vertex of the lattice $\mathbf{Z}^m$, such that neighbouring lines always intersect.
In this work, we present and study a novel class of discrete systems arising from billiards within confocal quadrics.
Unlike double reflection nets, these new systems are defined on a non-cubic lattice and they  assign hyperplanes of a given projective space to the verteces of the lattice.

This paper is organised as follows.
Section \ref{sec:drn} contains overview of confocal quadrics, related biliards, and double reflection nets.
In Section \ref{sec:hyper}, we introduce new discrete systems that are naturally related to the double reflection nets.
We derive their basic properties and discuss the cases when $m=2$ and $m=3$.

\section{Double reflection nets}
\label{sec:drn}
In this section, we give first a brief overview of most important notions regarding confocal families of quadrics and double reflection nets.
For more details on billiards within quadrics and confocal families see \cite{DragRadn2011book} and references therein.
On double reflection nets, see \cites{DragRadn2012jnmp,DragRadn2014bul,DragRadn2015umn}.

A \emph{pencil of quadrics} in the $d$-dimensional projective space $\mathbf{P}^d$ is a $1$-parameter family of quadrics:
\begin{equation}\label{eq:pencil}
\mathcal{Q}_{\lambda}\ :\ \big((A-\lambda I)x,x\big)=0,
\end{equation}
where $A$ is a symmetric non-degenerate operator and $I$ the identity operator.
A \emph{confocal family of quadrics} is a family dual to a pencil of quadrics:
\begin{equation}\label{eq:confocal}
\mathcal{Q}_{\lambda}^*\ :\ \big((A-\lambda I)^{-1}x,x\big)=0.
\end{equation}
When a confocal family of quadrics is given in the projective space, it is possible to define the billiard reflection off the quadrics from the family, see \cite{CCS1993}.
That definition is consistent with the billiard reflection in the Euclidean space.
By Chasles' theorem, each line in the space is touching exactly $d-1$ quadrics from (\ref{eq:confocal}) and these quadrics are preserved by reflections off quadrics from the confocal family.
Thus, a billiard trajectory within confocal quadrics always have $d-1$ \emph{caustics}.

Lines $\ell$, $\ell_1$, $\ell_2$, $\ell_{12}$ represent a \emph{double reflection configuration} if there are quadrics $\mathcal{Q}_{\alpha}$ and $\mathcal{Q}_{\beta}$ in (\ref{eq:confocal}) such that:
\begin{itemize}
\item
pairs $\ell$, $\ell_1$ and $\ell_2$, $\ell_{12}$ and satisfy reflection law off $\mathcal{Q}_{\alpha}^*$;
\item
pairs $\ell$, $\ell_2$ and $\ell_1$, $\ell_{12}$ and satisfy reflection law off $\mathcal{Q}_{\beta}^*$;
\item
four tangent planes at the reflection points are in a pencil.
\end{itemize}


A \emph{double reflection net} is a map $\varphi\ :\ \mathbf{Z}^m \to \mathcal{L}$, with $\mathcal{L}$ be the set of all lines in the projective space,
such that there exist $m$ quadrics $\mathcal{Q}_1^*$, \dots , $\mathcal{Q}_m^*$ from the confocal pencil, satisfying the following conditions:
\begin{itemize}
\item
the sequence $\{\varphi(\mathbf{n_0} + i\mathbf{e}_j)\}_{i\in\mathbf{Z}}$ represents a billiard trajectory within $\mathcal{Q}_j^*$,
 for each $j \in \{1,\dots,m\}$ and $\mathbf{n_0} \in\mathbf{Z}_m$;
 \item
the lines  $\varphi(\mathbf{n_0})$, $\varphi(\mathbf{n_0} + \mathbf{e}_i)$,  $\varphi(\mathbf{n_0} + \mathbf{e}_j)$,  $\varphi(\mathbf{n_0} + \mathbf{e}_i+\mathbf{e}_j)$ form a double reflection configuration, for all $i,j \in \{1,...,m\}$, $i \neq j$ and $\mathbf{n_0} \in \mathbf{Z}^m$.
\end{itemize}
In other words, for each edge in $\mathbf{Z}^m$ of direction $\mathbf{e}_i$, the lines corresponding to its vertices meet at $\mathcal{Q}_i^*$, while the four tangent planes at the intersection points, associated to an elementary quadrilateral, belong to a pencil.

The integrability of double reflection nets follows from the Six-pointed star theorem \cite{DragRadn2008}, which states that tere exist configurations consisting of twelve planes with the following properties:
\begin{itemize}
\item
The planes may be organized in eight triplets, such that each plane in a triplet is tangent to a different quadric from (\ref{eq:confocal}) and the three touching points are collinear.
Every plane in the configuration is a member of two triplets.

\item
The planes may be organized in six quadruplets, such that the planes in each quadruplet belong to a pencil and are tangent to two different quadrics from (\ref{eq:confocal}).
Every plane in the configuration is a member of two quadruplets.
\end{itemize}
Moreover, such a configuration is determined by three planes tangent to three different quadrics from (\ref{eq:confocal}), with collinear touching points.

\begin{figure}[h]
\centering
\psset{unit=1}
\begin{pspicture*}(-4,-3)(4,3)

\psset{viewpoint=-15 30 5,Decran=80,lightsrc=-5 10 50}

\psset{linewidth=0.01, linecolor=gray}

\psSolid[
object=line,
args=1 1 1 1 1 -1]
\psSolid[
object=line,
args=1 1 1 1 -1 1]
\psSolid[
object=line,
args=1 1 1 -1 1 1]
\psSolid[
object=line,
args=1 -1 -1 1 1 -1]
\psSolid[
object=line,
args=1 -1 -1 1 -1 1]
\psSolid[
object=line,
args=1 -1 -1 -1 -1 -1]
\psSolid[
object=line,
args=-1 -1 1 1 -1 1]
\psSolid[
object=line,
args=-1 -1 1 -1 1 1]
\psSolid[
object=line,
args=-1 -1 1 -1 -1 -1]
\psSolid[
object=line,
args=-1 1 -1 -1 1 1]
\psSolid[
object=line,
args=-1 1 -1 -1 -1 -1]
\psSolid[
object=line,
args=-1 1 -1 1 1 -1]

\psset{linewidth=0.05, linecolor=black}

\psSolid[object=line,
args=0 1 1 1 0 1]
\psSolid[object=line,
args=1 0 1 0 -1 1]
\psSolid[object=line,
args=0 -1 1 -1 0 1]
\psSolid[object=line,
args=-1 0 1 0 1 1]

\psSolid[object=line,
args=0 1 -1 1 0 -1]
\psSolid[object=line,linestyle=dashed,linewidth=0.02,
args=1 0 -1 0 -1 -1]
\psSolid[object=line,linestyle=dashed,linewidth=0.02,
args=0 -1 -1 -1 0 -1]
\psSolid[object=line,
args=-1 0 -1 0 1 -1]

\psSolid[object=line,
args=1 0 1 1 1 0 ]
\psSolid[object=line,
args=1 1 0 1 0 -1]
\psSolid[object=line,linestyle=dashed,linewidth=0.02,
args=1 0 -1 1 -1 0 ]
\psSolid[object=line,linestyle=dashed,linewidth=0.02,
args=1 -1 0 1 0 1]

\psSolid[object=line,
args=-1 0 1 -1 1 0 ]
\psSolid[object=line,
args=-1 1 0 -1 0 -1]
\psSolid[object=line,
args=-1 0 -1 -1 -1 0 ]
\psSolid[object=line,
args=-1 -1 0 -1 0 1]

\psSolid[object=line,
args=0 1 1 1 1 0 ]
\psSolid[object=line,
args=1 1 0 0 1 -1]
\psSolid[object=line,
args=0 1 -1 -1 1 0 ]
\psSolid[object=line,
args=-1 1 0 0 1 1]

\psset{linestyle=dashed,linewidth=0.02}

\psSolid[object=line,
args=0 -1 1 1 -1 0 ]
\psSolid[object=line,
args=1 -1 0 0 -1 -1]
\psSolid[object=line,
args=0 -1 -1 -1 -1 0 ]
\psSolid[object=line,
args=-1 -1 0 0 -1 1]

\psset{linestyle=solid,linewidth=0.05}

\psSolid[object=point, args=0 1 1]
\psSolid[object=point, args=0 -1 1]
\psSolid[object=point, args=0 1 -1]
\psSolid[object=point, args=0 -1 -1]
\psSolid[object=point, args=1 0 1]
\psSolid[object=point, args=1 0 -1]
\psSolid[object=point, args=-1 0 1]
\psSolid[object=point, args=-1 0 -1]
\psSolid[object=point, args=1 1 0]
\psSolid[object=point, args=-1 1 0]
\psSolid[object=point, args=1 -1 0]
\psSolid[object=point, args=-1 -1 0]

\end{pspicture*}
\caption{A configuration of planes from the Six-pointed star theorem.}\label{fig:cubo-oct}
\end{figure}
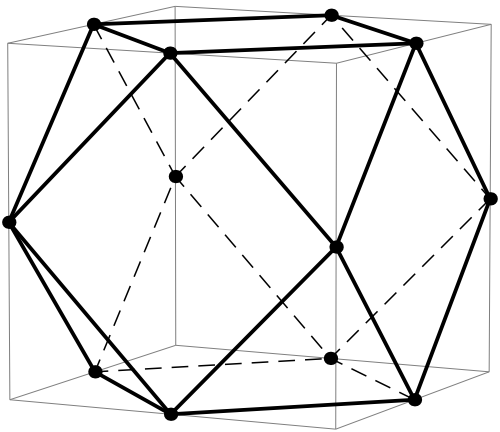

Such a configuration of planes is shown in Figure \ref{fig:cubo-oct}: each plane corresponds to a vertex of the cuboctahedron.
 These cuboctahedra will represent basic building blocks for lattices of hyper-planes that we introduce in the next section.

\section{Hyper-plane billiard nets}
\label{sec:hyper}
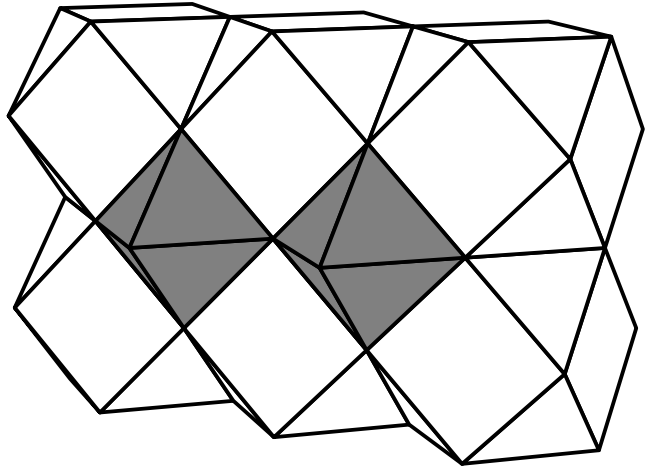
\begin{figure}[h]
\centering
\psset{unit=1}
\begin{pspicture*}(-4.5,-2)(4.5,4.5)

\psset{viewpoint=-15 30 8,Decran=30}

\psset{linewidth=0.05}

\psSolid[
object=cube,
a=3,
action=draw**,
trunccoeff=.5, trunc=all]
(3,0,0)

\psSolid[
object=cube,
a=3,
action=draw**,
trunccoeff=.5, trunc=all]
(3,0,3)

\psSolid[
object=cube,
a=3,
action=draw**,
trunccoeff=.5, trunc=all]

\psSolid[
object=cube,
a=3,
action=draw**,
trunccoeff=.5, trunc=all]
(0,0,3)

\psSolid[
object=cube,
a=3,
action=draw**,
trunccoeff=.5, trunc=all]
(-3,0,0)

\psSolid[
object=cube,
a=3,
action=draw**,
trunccoeff=.5, trunc=all]
(-3,0,3)

\psSolid[
object=octahedron,fillcolor=gray,
a=1.5,
action=draw**]
(1.5,1.5,1.5)

\psSolid[
object=octahedron,fillcolor=gray,
a=1.5,
action=draw**]
(-1.5,1.5,1.5)

\end{pspicture*}
\caption{Honeycomb consisting of cuboctahedra and octahedra.}\label{fig:honey3}
\end{figure}

Consider a lattice $\mathbf{Z}^m$ in $\mathbf{R}^m$.
That lattice generates a \emph{honeycomb}, that is a filling of the space by polytopes \cite{Coxeter}.
In this case, it is a regular honeycomb consisting of $m$-cubes.
The set of all midpoints of the edges of the cubes is
$$
\mathcal{M}^m=\bigcup_{1\le i\le m}\big(\mathbf{Z}^m+\frac12\mathbf{e}_i\big).
$$
The lattice $\mathcal{M}^m$ determines a honeycomb containing two types of polytopes (see \cite{Coxeter}):
\begin{itemize}
\item
\emph{rectified $m$-cubes} -- the verteces of each polytope of this kind are midpoints of the edges of an $m$-cube in the lattice $\mathbf{Z}^m$;
\item
\emph{cross polytopes} -- the verteces of each cross polytope are midpoints of all edges with a common endpoint in $\mathbf{Z}^m$.
\end{itemize}
This is an example of a convex uniform honeycomb \cite{Wells}.
For $m=3$, it is shown in Figure \ref{fig:honey3}, and for $m=2$ in Figure \ref{fig:tiling2}.
Integrable systems on such lattices were studied in \cites{KingS2003,KingS2006}

Assume that $\varphi \ :\ \mathbf{Z}^m\to\mathcal{L}$ is a given double reflection net.
For each $\mathbf{n_0}\in\mathbf{Z}^m$ and $i\in\{1,\dots,m\}$, the lines $\varphi(\mathbf{n_0})$ and $\varphi(\mathbf{n_0}+\mathbf{e}_i)$ are reflected to each other off the quadric $\mathcal{Q}_i^*$ from the confocal family.
We define the map
$$
\mathcal{H}\ :\ \mathcal{M}^m\to \mathbf{P}^{d*},
$$
such that it assigns to the midpoint of the edge $(\mathbf{n_0}, \mathbf{n_0}+\mathbf{e}_i)$ the tangent plane to $\mathcal{Q}_i^*$ at the point of reflection.
We introduce also the map
$$
\mathcal{P}\ :\ \mathcal{M}^m\to \mathbf{P}^{d*},
$$
which assigns the intersection point $\varphi(\mathbf{n_0})\cup\varphi(\mathbf{n_0}+\mathbf{e}_i)$  to the the midpoint of the edge $(\mathbf{n_0}, \mathbf{n_0}+\mathbf{e}_i)$.

Since each hyperplane of the space, except of a subset of measure $0$, is touching exactly one quadric from the given confocal family.
Thus, the touching point is uniquely determined and map $\mathcal{P}$ is uniquely determined by $\mathcal{H}$.
The inverse, to determine $\mathcal{H}$ when $\mathcal{P}$ is given, is not straightforward, since each point of the $d$-dimensiona space belongs to $d$ confocal quadrics.

From the construction, $\mathcal{H}$ and $\mathcal{P}$ have the following properties:
\begin{prop}
\begin{itemize}
\item
For each cross polytope of the honeycomb, $\mathcal{P}$ assigns to all its vertices collinear points.
Moreover, the points joined to the opposite vertices are on the same quadric from the confocal family.
\item
For each square $2$-face of any rectified $m$-cube, $\mathcal{H}$ assignes to its vertices hyperplanes that belong to one pencil and form a harmonic quadruple.
The hyperplanes corresponding to the opposite vertices are tangent to the same quadric from the confocal family.
\item
The hyperplanes assigned to any two adjacent verteces of a square $2$-face uniquely determine the hyperplanes assigned to the other two verteces.
\end{itemize}
\end{prop} 

\begin{rem}
For $m=2$, $\mathcal{M}^2$ determines a regular tesselation by squares, see Figure \ref{fig:tiling2}.
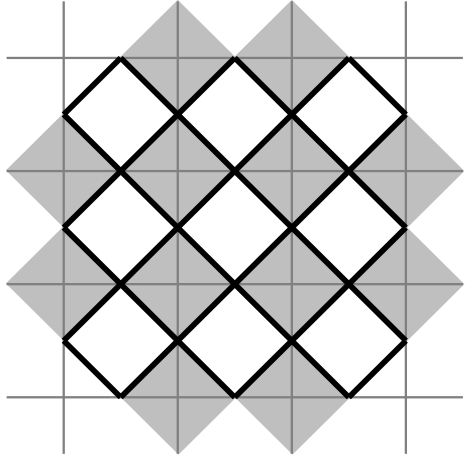
\begin{figure}[h]
\centering

\psset{unit=1.5}
\begin{pspicture}(-0.5,-0.5)(3.5,3.5)

\psset{fillstyle=solid, linecolor=gray!50, fillcolor=gray!50}
\psline(0,0.5)(.5,1)(0,1.5)(-0.5,1)

\psset{fillstyle=solid, linecolor=gray!50, fillcolor=gray!50}
\psline(0,1.5)(.5,2)(0,2.5)(-0.5,2)

\psset{fillstyle=solid, linecolor=gray!50, fillcolor=gray!50}
\psline(1,-0.5)(1.5,0)(1,.5)(0.5,0)

\psset{fillstyle=solid, linecolor=gray!50, fillcolor=gray!50}
\psline(1,0.5)(1.5,1)(1,1.5)(0.5,1)

\psset{fillstyle=solid, linecolor=gray!50, fillcolor=gray!50}
\psline(1,1.5)(1.5,2)(1,2.5)(0.5,2)

\psset{fillstyle=solid, linecolor=gray!50, fillcolor=gray!50}
\psline(1,2.5)(1.5,3)(1,3.5)(0.5,3)

\psset{fillstyle=solid, linecolor=gray!50, fillcolor=gray!50}
\psline(3,0.5)(3.5,1)(3,1.5)(2.5,1)

\psset{fillstyle=solid, linecolor=gray!50, fillcolor=gray!50}
\psline(3,1.5)(3.5,2)(3,2.5)(2.5,2)

\psset{fillstyle=solid, linecolor=gray!50, fillcolor=gray!50}
\psline(2,-0.5)(2.5,0)(2,0.5)(1.5,0)

\psset{fillstyle=solid, linecolor=gray!50, fillcolor=gray!50}
\psline(2,0.5)(2.5,1)(2,1.5)(1.5,1)

\psset{fillstyle=solid, linecolor=gray!50, fillcolor=gray!50}
\psline(2,1.5)(2.5,2)(2,2.5)(1.5,2)

\psset{fillstyle=solid, linecolor=gray!50, fillcolor=gray!50}
\psline(2,2.5)(2.5,3)(2,3.5)(1.5,3)

\psset{linecolor=gray,linewidth=0.02}

\psline(-0.5,0)(3.5,0)
\psline(-0.5,1)(3.5,1)
\psline(-0.5,2)(3.5,2)
\psline(-0.5,3)(3.5,3)

\psline(0,-0.5)(0,3.5)
\psline(1,-0.5)(1,3.5)
\psline(2,-0.5)(2,3.5)
\psline(3,-0.5)(3,3.5)

\psset{linecolor=black,linewidth=0.05}

\psline(0.5,0)(0,0.5)
\psline(1.5,0)(0,1.5)
\psline(2.5,0)(0,2.5)
\psline(3,0.5)(0.5,3)
\psline(3,1.5)(1.5,3)
\psline(3,2.5)(2.5,3)

\psline(0,0.5)(2.5,3)
\psline(0,1.5)(1.5,3)
\psline(0,2.5)(0.5,3)
\psline(0.5,0)(3,2.5)
\psline(1.5,0)(3,1.5)
\psline(2.5,0)(3,0.5)

\end{pspicture}
\caption{The lattice $\mathcal{M}^2$.}\label{fig:tiling2}
\end{figure}
For each square, if the values of $\mathcal{H}$ are given at two neighbouring vertices, it is possible to uniquely determine the values at the remaining two vertices.
The hyperplanes joined to the opposite points of a square are always tangent to the same quadric from the confocal pencil.
However, the discrete dynamics depends on the type of each square:
\begin{itemize}
\item
For the squares with vertices of the form 
$$
\mathbf{n_0}+\dfrac{\mathbf{e_1}}{2},
\quad
\mathbf{n_0}+\dfrac{\mathbf{e_2}}{2},
\quad
\mathbf{n_0}+\mathbf{e_1}+\dfrac{\mathbf{e_2}}2,
\quad
\mathbf{n_0}+\dfrac{\mathbf{e_1}}2+\mathbf{e_2},
$$
the corresponding hyperplanes are in a pencil and harmonically conjugated.
Such squares are white in Figure \ref{fig:tiling2}.
\item
For the squares with vertices of the form 
$$
\mathbf{n_0}+\dfrac{\mathbf{e_1}}{2},
\quad
\mathbf{n_0}+\dfrac{\mathbf{e_2}}{2},
\quad
\mathbf{n_0}-\dfrac{\mathbf{e_1}}{2},
\quad
\mathbf{n_0}-\dfrac{\mathbf{e_2}}2,
$$
the corresponding touching points are collinear, and in general not harmonically conjugated.
Such squares are gray coloured in Figure \ref{fig:tiling2}.
\end{itemize}
We denoted by $\mathbf{n_0}$ a point of $\mathbf{Z}^2$ and by $\mathbf{e_1}$, $\mathbf{e_2}$ unit coordinate vectors.
\end{rem}

It would be interested to write explicitely recursive relations that maps $\mathcal{H}$, $\mathcal{P}$ and double reflection nets satisfy on the lattice $\mathcal{M}^m$ and to see how they fit in the classification from \cite{ABS2003}.

\end{document}